\documentclass[a4paper,11pt]{article}
\usepackage{indentfirst,latexsym,bm,amsmath,amssymb,amsthm}
\usepackage[dvips]{graphicx}

\makeatletter\@addtoreset{equation}{section} \makeatother
\newtheorem{theorem}{Theorem}[section]
\newtheorem{lemma}[theorem]{Lemma}

 

\title{ Life-Span of Solutions to Critical Semilinear Wave Equations }
\author{Yi Zhou\thanks{  School of Mathematical Sciences, Fudan
University, Shanghai 200433, P. R. China;   Nonlinear Mathematical
Modeling and Methods Laboratory;    Shanghai Key Laboratory for
Contemporary Applied Mathematics  ({\tt Email: yizhou@fudan.ac.cn})
  } \and
Wei Han\thanks{School of Mathematical Sciences, Fudan University,
Shanghai 200433, P. R. China;     Department of Mathematics, North
University of China,  Taiyuan, Shanxi  030051, P. R. China  ({\tt
Email: sh\_hanweiwei1@126.com})}.}
\date{}
\begin{document}
\maketitle
\begin{abstract}
The final open part of  the famous Strauss conjecture on semilinear
wave equations of the form  $ \Box u=|u|^{p} $,  i.e.,    blow-up
theorem for the critical case in high dimensions  was solved by
Yordanov and Zhang \cite{Q. S. Zhang 2}, or Zhou  \cite{Y. Zhou 5}
independently. But the estimate for the lifespan, the maximal
existence time, of solutions was not clarified in both papers.
  Recently,  Takamura and
Wakasa \cite{Takamura and Wakasa}  have obtained the sharp upper
bound of the lifespan of the solution to the critical semilinear
wave equations,  and  their method is based on the method in
Yordanov and Zhang \cite{Q. S. Zhang 2}.   In this  paper, we give a
much simple proof of the result of   Takamura and Wakasa
\cite{Takamura and Wakasa}   by using the method in Y. Zhou \cite{Y.
Zhou 5} for space dimensions  $ n\geq 2$.

  Simultaneously,  this  estimate  of the life span  also  proves  the last
     open optimality  problem  of  the general theory for fully
    nonlinear wave equations with small initial data
 in the case $n=4$  and  quadratic nonlinearity(One can see Li and
 Chen\cite{Ta-tsien} for references on the whole history).

\par {\bf Keywords:} Semilinear wave equation;  Critical exponent;
Cauchy problem;  Blow up; Lifespan

\end{abstract}
\section{Introduction and Main Results}

In this paper, we will consider the blow up of solutions of the
  Cauchy  problems for the  following  semilinear wave
equations:

\begin{equation} \label{1.1}
\left \{
\begin{array}{lllll}
\Box u=|u|^{p},  \   \   \   \
(x,t)\in R^{n}\times (0,+\infty),  \  \      n\geq 2, \       \\
 t=0:  \    \         \
     u =\varepsilon f(x),  \   \    \     \      u_{t} =\varepsilon g(x),   \  \  \   x\in R^{n}, \\
\end{array} \right.
\end{equation}
where
  \begin{equation} \label{1.2}
   \Box=\partial_{t}^{2}-\sum\limits_{i=1}^{n}\partial_{x_{i}}^{2}
\end{equation}
  is the wave operator.
  and the initial values    $ f, g \in   C_{ 0} ^{\infty } ( R^{n}) $  satisfy
\begin{equation} \label{1.3}
   supp(f, \  g) \subset    \{  x | \      |x|\leq 1     \},
\end{equation}
$\varepsilon>0 $
   is a small parameter,   we assume that $ n\geq 2$,
 and $p=p_{0}(n)$  is the positive root of the quadratic
 equation  $$ (n-1)p^{2}-(n+1)p-2=0. $$
 The number  $p_{0}(n)$ is known as the critical exponent of problem
 \eqref{1.1}.   Since it divides $ (1, +\infty)$ into two subintervals
so that the following take place:  If $  p\in ( 1, p_{c}(n)]$,  then
  \eqref{1.1} has no global solution for
 nonnegative initial values; if
$  p\in ( p_{c}(n), +\infty)$,  then solutions with small( and
sufficiently regular ) initial values exist for all time ( see for
e.g. \cite{Strauss 2}   \cite{Q. S. Zhang 2}  and \cite{Takamura and
Wakasa}),  this is the famous  Strauss' conjecture, we only give a
brief summary here and refer the reader to  \cite{JZ1, Kato 1,
Ta-tsien 1, Lindblad, Rammaha, Takamura 1, Q. S. Zhang 1}.
   The case $n=3$ was first done by F. John
\cite{John 1} in 1979,  he showed that  when $n=3$ global solutions
always exist if $ p>p_{0}(3)=1+\sqrt{2} $ and initial data are
suitably small, and moreover, the global solutions do not exist if
$1<p<p_{0}(3)=1+\sqrt{2} $ for any nontrivial choice of $f$ and $g$.
 The number $p_{0}(3)=1+\sqrt{2} $ appears to have first arisen in
Strauss' work on low energy scattering for the nonlinear
Klein-Gordon equation \cite{Strauss 1}. This led him to conjecture
that when $n\geq 2$ global solutions of \eqref{1.1} should always
exist if initial data are sufficiently small and $p$ is greater than
a critical power $p_{0}(n)$. The conjecture was verified when $n=2$
by R. T. Glassey \cite{Glassey 2}.   In higher space dimensions, the
case $n=4$  was proved by Y. Zhou \cite{Y. Zhou 3} and Lindblad and
Sogge \cite{LS} proved the case $ 3\leq n \leq 8 $.   Later V.
Georgiev, H. Lindblad and C. Sogge \cite{Georgiev} showed that when
$n\geq 4$ and $ p_{0}(n)<p\leq\frac{n+3}{n-1}$, \eqref{1.1} has
global solutions for small initial values (see also \cite{LS} and
\cite{Tataru}). Later,  a simple proof was given by Tataru
\cite{Tataru} in the case $p>p_{0}(n) $ and $n\geq 4$. R. T. Glassey
\cite{Glassey 1} and  T. C. Sideris \cite{Sideris} showed the
blow-up result of $1<p<p_{0}(n)$ for $n=2 $  and  all $n\geq 4 $,
respectively.
     On the other hand,  for
the critical case $p=p_{0}(n)$, it was shown by Schaeffer
\cite{Schaeffer} that the critical power also belongs to the blowup
case  for small data when $n=2,  3$  (see also \cite{Takamura, Y.
Zhou 1,Y. Zhou 2}).
  When $ n\geq 4$,  the blow-up problem
  for the critical wave equations problem
  was solved by Yordanov and Zhang
\cite{Q. S. Zhang 2} and Zhou \cite{Y. Zhou 5} by different methods
  respectively.   But the sharp estimate for the lifespan, the maximal
existence time, of solutions was not clarified in both papers.

 Recently,  Takamura and
Wakasa \cite{Takamura and Wakasa}  have obtained the sharp upper
bound  $ T(\varepsilon) \leq
              \exp ({ B\varepsilon^{-p(p-1) } })  $ of the lifespan of the solution to the critical semilinear
wave equations for $ n\geq 4 $,   and method is based on the method
in Yordanov and Zhang \cite{Q. S. Zhang 2}.  For the sharp upper
bound in the low dimensional csae,   $ n=2 $ was obtained by   Zhou
\cite{Y. Zhou 2}
 and $ n=3 $ was obtained by Zhou \cite{Y. Zhou 1} much earlier,
 and
 Takamura \cite{Takamura} gave an another uniform proof for the
 case  $n=2, 3 $.
 In this paper,  we give a much simple
proof of their result based on  the method in Y. Zhou \cite{Y. Zhou
5} for space dimensions $n \geq2$ .

 Simultaneously,   the lifespan $T(\varepsilon)$ of solutions of
  $  \Box u =  u^{2} $ in  $ R^{4 } \times [ 0, \infty)
  $ with the initial data  $  u(0, x)= \varepsilon f(x),
   \   u_{ t}(0, x)= \varepsilon g(x)$ of a small parameter $ \varepsilon >0
   $, compactly supported smooth functions $f$ and $g$, has an upper
   bound estimate
   $ T(\varepsilon) \leq
              \exp ({ B\varepsilon^{-2 } }) $ with a positive constant $ B $ independent
              of $\varepsilon
  $, which belongs to the same kind of the lower bound of the
  lifespan.     This upper bound  proves  the last
              open sharpness  problem
    on
the lower bound
 for the lifespan of solutions to  fully
              nonlinear wave equations   with small initial data
 in the case of $n=4$  and  quadratic nonlinearity(One can see Li Ta-Tsien and Chen Yunmei
 \cite{Ta-tsien}  for references on the whole history).

   We
define ``life span" $T(\varepsilon)$ of the solutions of \eqref{1.1}
   to be the largest value such that solutions exist for
$x\in R^{n}$,  $0\leq t< T(\varepsilon)$.

 For problem  \eqref{1.1},  we consider compactly supported nonnegative data $(f,g)\in
  C_{ 0} ^{\infty } ( R^{n})$,  $n\geq 2$ and satisfy
\begin{equation}\label{1.4}
f(x)\geq 0,  \   \    g(x)\geq 0, \    \forall  x \in  R^{ n},  \  \
  \      \mbox { and}  \  \    g(x)\not\equiv 0.
\end{equation}
  We establish the following theorem for \eqref{1.1}:

\noindent\begin{theorem}\label{thm:1.1}   For Cauchy problem
\eqref{1.1},   let  initial values  $f$ and $g$ satisfies
\eqref{1.3} and  \eqref{1.4},  space dimensions  $ n\geq 2 $, $
p=p_{0}(n)$,    suppose that Cauchy problem \eqref{1.1}
  has a solution  satisfying that
   \begin{equation} \label{1.5}
    u \in
C([0,T), \ H^{1}(R^{n}) )  \cap   C^{ 1}([0,T), \ L^{2}(R^{n}) )
\cap  C([0,T), \ L^{p}(R^{n}) ),
 \end{equation}
 \begin{equation} \label{1.6}
  supp u   \subset    \{ ( x, t)  | \      |x|\leq t+ 1     \}.
 \end{equation}
Then  lifespan
  $T<\infty$, and there exists a
 positive constant $B $ which is independent of $ \varepsilon$
 such that
\begin{equation}
\begin{array}{ll}
  T(\varepsilon) \leq
              \exp ({ B\varepsilon^{-p(p-1) } }).
\end{array}
\end{equation}
Especially, when  $ n=4, $   $ p=p_{0}(4)=2 $,  corresponding to the
case of $ n=4$  and  quadratic nonlinearity(see Li Ta-Tsien and Chen
Yunmei
 \cite{Ta-tsien}), so we prove that   $ T(\varepsilon)
\leq  \exp ({ B\varepsilon^{-2} }) $.
\end{theorem}

The rest of the paper is arranged as follows. We state several
preliminary propositions in Section 2, Section 3 is devoted to the
proof  of our Theorem \ref{thm:1.1}.

\section{Preliminaries} \vskip .5cm

 In this section, we consider Cauchy problem \eqref{1.1},
  where the initial data  satisfy \eqref{1.3}, and we consider  the case  $n\geq 2
  $,  $ p=p_{ c}(n) $.

To prove the main results in this paper, we will employ  some
important lemmas,

 \noindent
\begin{lemma} \label{thm:2.1}  Suppose that   $K (t) $  and $ h(t)$
are all positive $ C^{ 2}$ functions  and satisfies

  \begin{equation} \label{2.1}
   a(t)K''(t) + K'(t) \geq  b(t) K^{ 1+\alpha } (t),   \    \      \
   \forall  t \geq 0,
 \end{equation}

 \begin{equation} \label{2.2}
   a(t)h''(t) + h'(t) \leq   b(t) h^{ 1+\alpha } (t),   \    \      \
   \forall  t \geq 0,
 \end{equation}
 where  $ \alpha \geq 0 $, and
 \begin{equation} \label{2.3}
   a(t),   \    b(t) > 0, \    \    \        \forall  t\geq  0.
   \end{equation}
 If
 \begin{equation} \label{2.4}
  K(0) >h(0),  \   \     \          K'(0)  \geq  h'(0).
 \end{equation}
 Then we have
\begin{equation} \label{2.5}
  K'(t)>   h'(t),  \   \    \    \forall  t>  0,   \end{equation}
 which imply
\begin{equation} \label{2.6}
 K(t)>   h(t),  \   \    \    \forall  t\geq  0.
 \end{equation}
\end{lemma}

\noindent\begin{proof}\
 Without loss of generality, we can assume that
 $$  K'(0)>h' (0).          $$  Otherwise, if $ K'(0)= h' (0)  $,
 then we can obtain from    \eqref{2.1} and  \eqref{2.2},     $ K''(0)>h'' (0)
 $,  so there exists a positive constant  $ \delta_{ 0} $ such that
   $$  K'(t)>h' (t), \   \    \      \forall  0<t \leq  \delta_{ 0},    $$
  so we only need to take  $ \delta_{ 0}  $ as the initial time.

  By the continuity,  assume \eqref{2.5} is not true, then there
  exists a positive constant   $  t^{* }>0, $  such that

$$
  \left\{
     \begin{array}{ll}
      K'( t ) > h'(t),   \   \    \         \forall   0\leq  t < t^{
      *},                \cr\noalign{\vskip4mm}
         K'( t^{ *} ) = h'(t^{ *}).
\end{array}
   \right.
$$
 So we obtain that
 $$  K''( t^{ *} ) \leq  h''(t^{ *}) .   $$
 On the other hand, by \eqref{2.4}, we have  $  K( t^{ *} )   >  h(t^{ *}) ,    $
  so by  \eqref{2.1} and  \eqref{2.2},  we can get
  $$    K''( t^{ *} ) >  h''(t^{ *}).      $$
 This is a contradiction, so \eqref{2.5} holds, furthermore,
 \eqref{2.6} also holds.
\end{proof}
 \noindent
\begin{lemma} \label{thm:2.2}  Suppose that  Cauchy problem \eqref{1.1}  has a solution $ u$, such that
    \begin{equation} \label{2.7}
    u \in
C([0,T), \ H^{1}(R^{n}) )  \cap   C^{ 1}([0,T), \ L^{2}(R^{n}) )
\cap  C([0,T), \ L^{p}(R^{n}) ),
 \end{equation}
 \begin{equation} \label{2.8}
  supp u   \subset    \{ ( x, t)  | \      |x|\leq t+ 1     \}.
 \end{equation}
If the initial data satisfies that
 \begin{equation} \label{2.9}
   \displaystyle\int_{ R^{ n}}  \phi_{ 1} ( x) f(x) dx\geq  0,   \
    \      \         \displaystyle\int_{ R^{ n}}  \phi_{ 1} ( x) g(x) dx\geq  0,
 \end{equation}
 and they are not identically zero.  where  $ \phi_{ 1} ( x) $ is
 defined as
  $$  \phi_{1}(x)=  \displaystyle\int_{S^{n-1}} e^{x\cdot \omega} d\omega.  $$
 Then we have
 $$   \displaystyle\int_{ R^{ n}}   | u( t, x)|^{p } dx  \geq   C_{ 1} \varepsilon^{ p} ( 1+t )^{ n-1-\frac{ n-1}{ 2} p  } ,
    \   \      t\geq 1.    $$
\end{lemma}

\noindent\begin{proof}\   This is exactly ( 2.5' ) in Yordanov and
Zhang\cite{Q. S. Zhang 2}, see also  \cite{Q. S. Zhang 1}.
\end{proof}

For the convenience of describing  the following
  lemma, let us introduce  a special function,  we seek a solution
of the linear wave equation
 \begin{equation} \label{2.10}
   \Box \phi=0
\end{equation}
 on the domain $  | x|\leq t, \    t\geq 0 $ of the  following  form
 \begin{equation} \label{2.11}
     \phi= \phi_{ q}= ( t+|x| )^{ -q} h \left( \frac{ 2| x| }{ t+|x |  }
      \right),
\end{equation}
 where  $ q>0 $,  we substituting \eqref{2.11} in \eqref{2.10},  by
 easy computation one obtain that $ h=h_{ q} $ satisfies the
 ordinary  differential equation
 \begin{equation} \label{2.12}
    z(1-z) h''(z)  +  \left[ n-1-\left( q+ \frac{ n+1 }{ 2}  \right) z \right] h'(z) - \frac{n-1  }{ 2} q h(z)
    =0,
\end{equation}
 where $ n $ stands for the space dimensions.  Therefore, we can
 take
 \begin{equation} \label{2.13}
    h_{ q}(z) = F( q, \frac{ n-1  }{2 },  n-1, z  ),
\end{equation}
where $F$ is the hypergeometric function defined by
$$  F( \alpha,  \beta,  \gamma,  z   )= \sum\limits_{ k=0 } ^{ \infty }
\frac{  ( \alpha )_{ k} ( \beta )_{ k}   }{  k! ( \gamma )_{ k }  }
z^{ k},  \     \      \             | z| <1,
$$
 with   $   ( \lambda )_{ 0} =1,  \   ( \lambda )_{ k}= \lambda ( \lambda +1 )\cdots ( \lambda +k-1 ),   k\geq 1.
 $   For $ \gamma > \beta >0  $,  we have the formula

  $$   F( \alpha,  \beta,  \gamma,  z )= \frac{ \Gamma ( \gamma)  }{ \Gamma( \beta )  \Gamma ( \gamma-\beta ) }
         \displaystyle\int_{0}^{1 } t^{ \beta-1 } ( 1-t )^{  \gamma-\beta-1 } ( 1-zt )^{ -\alpha } dt,    \   | z |<1. $$
Thus
\begin{equation} \label{2.14}
    h_{ q}(z) = \frac{ \Gamma ( n-1 )       }{  \Gamma^{ 2}( \frac{n-1}{ 2 }) }  \displaystyle\int_{0}^{1 }
       t^{ \frac{ n-3 }{ 2} } ( 1-t )^{ \frac{ n-3 }{ 2} } ( 1-zt )^{ -q}
       dt.
\end{equation}
 So we have
 \begin{equation} \label{2.15}
    h (z) >0,   \    \       0\leq  z< 1.
\end{equation}
 Moreover,  when
 \begin{equation} \label{2.16}
         0< q < \frac{ n-1 }{ 2} ,
\end{equation}
  $ h(z) $  is continuous at  $ z=1$.   Thus
 \begin{equation} \label{2.17}
          C_{0}^{ -1}\leq  h(z) \leq  C_{ 0},  \   \         0 \leq
          z\leq 1,
\end{equation}
  where $  C_{0} $ is a positive constant.   When
 $$
          q > \frac{ n-1 }{ 2} ,
    $$
 $ h(z) $  behaves like  $  ( 1-z)^{ -( q-\frac{ n-1 }{ 2} ) }  $
 when $ z $ is close to $ z=1$.   Thus

 \begin{equation} \label{2.18}
            C_{0}^{ -1} ( 1-z )^{  \frac{ n-1 }{ 2}-q } \leq  h(z)
            \leq  C_{ 0} ( 1-z )^{  \frac{ n-1 }{ 2}-q },
\end{equation}
for some positive constant $ C_{0} $.   Furthermore, One can easily
verify
 \begin{equation} \label{2.19}
    \partial _{ t} \phi_{ q}(t, x)= -q  \phi_{ q+1 }(t, x).
\end{equation}
 Based on these facts,   we have the following

\noindent \begin{lemma}\label{thm:2.3}   Consider Cauchy problem
\eqref{1.1},  where  initial data satisfies \eqref{1.3} and
\eqref{1.4},  space dimensions  $ n\geq 2 $,   $ p=p_{ 0}(n) $, and
its solution $ u $ also satisfies  \eqref{1.6}.   Let
 \begin{equation} \label{2.20}
     G(t) =   \displaystyle\int_{0}^{t } ( t-\tau ) ( 1+\tau )  \displaystyle\int_{ R^{
     n}} \Phi_{ q}( \tau, x ) | u( \tau, x)|^{p } dxd\tau ,
\end{equation}
 where   $ q= \frac{ n-1 }{ 2}-\frac{ 1}{ p}$,   $    \Phi_{ q}( \tau, x ) = \phi_{ q}( \tau+2, x ).
 $   Then we have
 \begin{equation} \label{2.21}
     G'(t) >  K_{0}  \left(  \ln (2+t )  \right)^{ -(p-1)} (2+t ) \left(  \displaystyle\int_{0}^{t } (2+ \tau)^{ -3} G( \tau
     ) d\tau   \right)^{ p},   \     \         t\geq 1,
\end{equation}
  where  $ K_{ 0} $ is a constant which is independent of $ \varepsilon
  $.
\end{lemma}

\noindent\begin{proof}\  By a simple calculation, we obtain
 \begin{equation} \label{2.22}
     G'(t) =   \displaystyle\int_{0}^{t }  ( 1+\tau )  \displaystyle\int_{ R^{
     n}} \Phi_{ q}( \tau, x ) | u( \tau, x)|^{p } dxd\tau ,
\end{equation}
 \begin{equation} \label{2.23}
     G''(t) =   (1+t)  \displaystyle\int_{ R^{
     n}} \Phi_{ q}( t, x ) | u( t, x)|^{p } dx,
\end{equation}
 multiplying the equation  \eqref{1.1}  by  $  \Phi_{ q}( t, x ) $,   $  q=\frac{ n-1 }{ 2}-\frac{ 1}{ p},
 $  and  using integration by parts, we have
   $$    \displaystyle\int  \Phi_{ q} ( u_{ tt} -\Delta u   ) dx =  \displaystyle\int  \Phi_{
   q} | u|^{ p} dx,   $$
 since
 $$    \displaystyle\int  \Phi_{ q} \Delta u  dx =  \displaystyle\int  \Delta \Phi_{
   q} u   dx =   \displaystyle\int   \Phi_{
   q tt} u   dx,   $$
so we have
  $$    \displaystyle\int  \Phi_{ q} ( u_{ tt} -\Delta u   ) dx =  \displaystyle\int  \left(   \Phi_{
   q} u_{ tt}  -   \Phi_{  q tt} u \right)   dx =  \displaystyle \frac{ d }{ dt } \displaystyle\int
      \left(  \Phi_{ q} u_{ t} -  \Phi_{ qt}u \right)  dx ,   $$
 and by \eqref{2.19}, we have
$$     \displaystyle\int  \left(  \Phi_{ q} u_{ t} -  \Phi_{ qt}u \right)  dx =
   \displaystyle \frac{ d }{ dt }  \displaystyle\int    \Phi_{ q} u dx- 2  \displaystyle\int   \Phi_{ qt}u dx
    =   \displaystyle \frac{ d }{ dt }  \displaystyle\int    \Phi_{ q} u dx +2q  \displaystyle\int  \Phi_{ q+1} u dx,    $$
 therefore we have

$$   \displaystyle \frac{ d^{ 2} }{ dt^{ 2} } \displaystyle\int   \Phi_{
q} u dx+   2q \displaystyle \frac{ d }{ dt }  \displaystyle\int
\Phi_{ q +1} u dx =  \displaystyle\int  \Phi_{ q} | u|^{ p} dx,  $$
 Integrating  the above expression  three times about $t$  from $ 0$  to $ t$,  we obtain
$$
\begin{array}{ll}
 & \displaystyle\int_{0}^{t } \displaystyle\int_{ R^{ n}} \Phi_{
q} u dx d\tau + 2 q \displaystyle\int_{0}^{t }  ( t-\tau )
\displaystyle\int_{ R^{ n}}   \Phi_{ q +1} u dxd\tau
  \cr\noalign{\vskip2mm}
  &  = \displaystyle \frac{ 1 }{2 } \displaystyle\int_{0}^{t } ( t-\tau )^{ 2}  \displaystyle\int_{ R^{
 n}}   \Phi_{ q}  | u|^{ p} dx d\tau +   \varepsilon t \displaystyle\int_{
 R^{ n}}  \Phi_{ q} ( 0,x ) f(x) dx +  \frac{\varepsilon t^{ 2}}{ 2
 } \displaystyle\int_{
 R^{ n}}    \left(  \Phi_{ q} ( 0,x ) g(x) + q  \Phi_{ q +1} ( 0,x )f(x)
 \right) dx.
 \end{array}
    $$
 Based on the positivity of the  $ \Phi_{ q}$ and the initial
 values,  we have
 \begin{equation} \label{2.24}
      \displaystyle\int_{0}^{t } \displaystyle\int_{ R^{ n}}  \Phi_{
q} u dx d\tau  +  2 q \displaystyle\int_{0}^{t }  ( t-\tau )
\displaystyle\int_{ R^{ n}}   \Phi_{ q +1} u dxd\tau > \displaystyle
\frac{ 1}{2} \displaystyle\int_{0}^{t } ( t-\tau )^{ 2 }
\displaystyle\int_{ R^{ n}} \Phi_{ q} | u|^{ p} dx d\tau,
\end{equation}
by using the finite propagation speed of  waves \eqref{2.8}, H$
\ddot{o} $lder's inequality , and noting \eqref{2.22}, we have
 $$   \displaystyle\int_{0}^{t } \displaystyle\int_{ R^{ n}}  \Phi_{
 q} u dxd\tau  \leq    \left( G'(t)  \right)^{\frac{ 1}{ p} }
 \left(  \displaystyle\int_{0}^{t } \displaystyle\int_{  |x|\leq \tau +1  }
  \Phi_{ q} ( 1+ \tau)^{  -\frac{p' }{ p} } dxd\tau   \right)^{ \frac{ 1}{ p'}
  },      $$
 where  $  \frac{ 1}{ p} +  \frac{ 1}{ p'}=1$.
 Since   $  0<q<  \frac{ n-1}{ 2}$,  so we have    $  \Phi_{ q} \sim
(1+ \tau )^{ -q} $,  therefore, we have
$$   \displaystyle\int_{0}^{t } \displaystyle\int_{  |x|\leq \tau +1
}  \Phi_{ q} ( 1+ \tau)^{  -\frac{p' }{ p} } dxd\tau  \leq   c
\displaystyle\int_{0}^{t } ( 1+ \tau )^{ -q+ n- \frac{p' }{ p} }
d\tau,     $$     since  $  p=p_{0}(n), $   $ q= n-1-\frac{ 2 }{
p-1}$, so we have
  $$   n-q-  \frac{p' }{ p}=1+ \frac{p' }{ p}.        $$
   Thus we obtain,  when  $ t\geq 2 $,
 \begin{equation} \label{2.25}
    \displaystyle\int_{0}^{t } \displaystyle\int_{ R^{ n}}  \Phi_{ q} u dxd\tau \leq
     C   \left( G'(t)  \right)^{\frac{ 1}{ p} } ( 1+t)^{ \frac{ 2  }{ p' } + \frac{ 1 }{ p} } =
      C   \left( G'(t)  \right)^{\frac{ 1}{ p} } ( 1+t)^{  1 + \frac{ 1 }{ p'} },
 \end{equation}
in same way, because   $  q+1>  \frac{ n-1 }{ 2} $,   so we have
 $  \Phi_{ q+1 } \sim
(2+ \tau )^{ -\frac{ n-1}{ 2}} ( 2+ \tau -|x | )^{ -( q+1- \frac{
n-1}{ 2} )   } $.
 So by H$ \ddot{o} $lder's  inequality, we can get
$$   \displaystyle\int_{0}^{t }  ( t-\tau ) \displaystyle\int_{ R^{
n}}   \Phi_{ q+1} u dxd\tau \leq   \left( G'(t)  \right)^{\frac{ 1}{
p} }  \left( \displaystyle\int_{0}^{t } ( t-\tau )^{ p'}
\displaystyle\int_{  |x|\leq \tau +1 } \Phi_{ q}  \left( \frac{
\Phi_{ q +1} }{
 \Phi_{ q} }
    \right)^{ p'} ( 1+\tau )^{-\frac{p' }{ p}   }  dxd\tau
   \right)^{ \frac{ 1}{ p' }  },
 $$
 and we know that
$$  \displaystyle\int_{  |x|\leq \tau +1 } \Phi_{ q}  \left( \frac{
\Phi_{ q +1} }{   \Phi_{ q} }  \right)^{ p'} ( 1+\tau )^{-\frac{p'
}{ p}   }  dx \leq  C (1+ \tau)^{ n-1+q( p'-1) -\frac{ n-1 }{2 }p'
-\frac{p' }{ p} } \displaystyle\int_{0}^{ 1+t }  ( 2+t-r )^{ -p'(
q+1- \frac{ n-1 }{2 }  ) } dr,    $$ we can easily verify that
 $$  p' (q+1- \frac{ n-1 }{2 }   )=1,    $$
  $$ n-1+q( p'-1 ) - \frac{ n-1 }{2 }p'- \frac{ p' }{ p}=0,     $$
so we have
  \begin{equation} \label{2.26}
   \displaystyle\int_{0}^{t }  ( t-\tau ) \displaystyle\int_{ R^{
n}}   \Phi_{ q+1} u dxd\tau \leq
  C  \left( G'(t)  \right)^{\frac{ 1}{
p} }  ( 1+t )^{ \frac{ 1+p' }{ p' }  }  \left(  \ln ( 2+t) \right)^{
\frac{ 1}{ p'}  }.
   \end{equation}
 Since when  $ t\geq 1$,     $  \ln ( 2+t)>1$.  and
  by   \eqref{2.24},   \eqref{2.25} and  \eqref{2.26},  and using
 the expression  \eqref{2.23} of $ G''(t) $,
  we get
$$  \left( G'(t)  \right)^{\frac{ 1}{ p} }  \left(  \ln ( 2+t) \right)^{
\frac{ 1}{ p'}  } ( 1+t )^{ 1+ \frac{ 1}{ p'} } \geq
  C \displaystyle\int_{0}^{t }  ( t-\tau )^{ 2} ( 1+\tau  )^{ -1} G''( \tau
  ) d\tau,
$$
 and we use integration by parts twice,  we can get
$$   \displaystyle\int_{0}^{t }  ( t-\tau ) ^{ 2} ( 1+\tau  )^{ -1}
 G''( \tau  ) d\tau =   \displaystyle\int_{0}^{t } \partial_{ \tau}^{
 2} \left[ ( t-\tau )^{ 2} ( 1+\tau  )^{ -1}   \right] G( \tau )
 d\tau,  $$
 by simple calculation, we have
  $$   \partial_{ \tau}^{
 2} \left[ ( t-\tau )^{ 2} ( 1+\tau  )^{ -1}   \right] = 2 ( 1+\tau  )^{
 -3} (t+1 )^{ 2}. $$
 So we have
$$  \left( G'(t)  \right)^{\frac{ 1}{ p} }  \left(  \ln ( 2+t) \right)^{
\frac{ 1}{ p'}  } ( 1+t )^{ 1+ \frac{ 1}{ p'} } \geq
   C( t+1)^{ 2}  \displaystyle\int_{0}^{t } ( 1+\tau  )^{
 -3} G( \tau ) d\tau ,    \    \        t\geq  1. $$
 This imply that
 $$   G'(t) \geq   C  \left(  \ln ( 2+t) \right)^{
   -( p-1 )  } ( 1+t ) \left( \displaystyle\int_{0}^{t } ( 1+\tau  )^{
 -3} G( \tau ) d\tau       \right)^{ p},  \   \       t\geq  1.       $$
Thus the proof of Lemma \ref{thm:2.3} is complete.
\end{proof}

\section{The proof of Theorem  \ref{thm:1.1}} \vskip .5cm

By Lemma \ref{thm:2.3},   we can get  \eqref{2.21}.  Let

 \begin{equation} \label{3.1}
     H(t) =   \displaystyle\int_{0}^{t }  ( 2+\tau )^{ -3 } G( \tau
     ) d\tau,
\end{equation}
 then
 \begin{equation} \label{3.2}
     H'(t) =    ( 2+ t )^{ -3 } G( t
     ),
\end{equation}
therefore we have
 \begin{equation} \label{3.3}
     G(t) =    ( 2+ t )^{ 3 } H'(t).
\end{equation}
   So  by \eqref{2.21},   we  obtain
 \begin{equation} \label{3.4}
         \left(  ( 2+ t )^{ 3 } H'(t) \right)'>  K_{ 0}  \left(  \ln ( 2+t) \right)^{
   -( p-1 )  }  ( 2+t ) H^{p }(t).
\end{equation}
By Lemma \ref{thm:2.2},  and    $  \Phi_{ q} \sim (1+ \tau )^{ -q}
$,     we can obtain
 \begin{equation} \label{3.5}
   \begin{array}{ll}
     G(t)  & =   \displaystyle\int_{0}^{t } ( t-\tau ) ( 1+\tau )  \displaystyle\int_{ R^{
     n}} \Phi_{ q}( \tau, x ) | u( \tau, x)|^{p } dxd\tau    \cr\noalign{\vskip3mm} &
   \geq   c  \displaystyle\int_{0}^{t } ( t-\tau ) ( 1+\tau )^{ 1-q
   }  \displaystyle\int_{ R^{
     n}}  | u( \tau, x)|^{p } dxd\tau  \cr\noalign{\vskip3mm} &
     \geq
      c \varepsilon^{p } \displaystyle\int_{0}^{t } ( t-\tau )  ( 1+\tau )^{
      1-q +n-1-  \frac{ n-1}{ 2} p
   }  d\tau  ,
      \end{array}
\end{equation}
and since   $ p= p_{ 0}(n)$,  we can easily check that
   $$ 1-q +n-1-  \frac{ n-1}{ 2} p=0,    $$
 so we have
   $$  G(t) \geq  C  \varepsilon^{p }  t^{ 2}.   $$
 Therefore
 \begin{equation} \label{3.6}
     H(t) \geq  c  \varepsilon^{p }  \displaystyle\int_{0}^{t } ( 2+ \tau )^{ -3 }  \tau^{ 2} d\tau
      \geq  c_{ 0} \varepsilon^{p } \ln ( 2+t ),   \     \     \
          t\geq  2,
\end{equation}
and
 \begin{equation} \label{3.7}
     H'(t) \geq  c  \varepsilon^{p }  ( 2+ t )^{ -3 } t^{ 2}
      \geq  c_{ 0} \varepsilon^{p } ( 2+t )^{ -1},   \     \     \
          t\geq  2.
\end{equation}
 We open the expression   \eqref{3.4},   we can get
 $$           ( 2+ t )^{ 2 } H''(t) +   3 ( 2+t )H'(t)
 >  K_{ 0}  \left(  \ln ( 2+t) \right)^{
   -( p-1 )  }   H^{p }(t).
   $$
Let us make a transformation   $   t+2=  \exp ( \tau ) $,    and
define $ H_{ 0}( \tau  ) =  H( \exp ( \tau )-2   ) =H(t)$,  one has
$$  H_{ 0}'( \tau )= H'(t) \frac{dt }{ d\tau } = ( t+2)H'(t),       $$
$$  H_{ 0}''( \tau )=    \left( (t+2)H'(t) \right)' (t+2)= ( t+2 )^{ 2} H''(t) +  ( t+2 )H'(t).   $$
  By  \eqref{3.4},   we have
$$
  \left\{
     \begin{array}{ll}
     H_{0}''( \tau ) +2 H_{0}'( \tau ) >  K_{ 0} \tau ^{ -( p-1 )}
     H_{0}^{ p} ( \tau ),   \    \     \cr\noalign{\vskip3mm}
     H_{0}( \tau ) \geq    C_{ 0} \varepsilon^{p } \tau ,     \    \     \cr\noalign{\vskip 3mm}
         H_{0}'( \tau ) \geq  C_{ 0}  \varepsilon^{p } .
\end{array}
   \right.
$$
Let    $$  H_{ 1}(s) =  \varepsilon^{p^{2} -2p } H_{0} (
\varepsilon^{   -p( p-1 ) } s).
$$
 Then we have

$$
  \left\{
     \begin{array}{ll}
     \varepsilon^{   p( p-1 ) }  H_{1}''( s ) +2 H_{1}'( s ) >  K_{ 0} s ^{ -( p-1 )}
     H_{1}^{ p} ( s ),   \    \     \cr\noalign{\vskip 3mm}
     H_{1}( s ) \geq    C_{ 0} s ,     \    \     \cr\noalign{\vskip 3mm}
         H_{1}'( s ) \geq  C_{ 0}.
\end{array}
   \right.
$$

We take  $ s_{ 0} $,  $ \delta $ independent of  $ \varepsilon $,
and   $K_{ 0} ,   C_{ 0} \ll  s_{ 0} \ll  \frac{ 1  }{ \delta } $.
Let       $ H_{ 2}(s) = s H_{ 3}(s)$,    and    $ H_{ 3}(s) $
satisfies
 \begin{equation} \label{3.8}
  \left\{
     \begin{array}{ll}
      H_{3}'( s ) =  \delta H_{3}^{ \frac{ p+1 }{ 2}  }( s ),     \    \    s\geq s_{ 0},    \cr\noalign{\vskip3mm}
         H_{3}( s_{ 0} ) =  C_{ 0}/4.
\end{array}
   \right.
\end{equation}
Then we have
$$  H_{2}'( s )= H_{3}(s) +s H_{3}'(s) = H_{3}(s)+ \delta s H_{3}^{ \frac{ p+1 }{ 2} }( s),    $$
$$  H_{2}''( s )= 2 \delta  H_{3}^{ \frac{ p+1 }{ 2} }(s)  + \delta^{ 2}s \left( \frac{ p+1 }{ 2}  \right) H_{3}^{ p}(s),    $$
 so we have
$$     \varepsilon^{   p( p-1 ) }  H_{2}''( s ) + 2 H_{2}'( s ) =  \varepsilon^{   p( p-1 )
}  \left( \frac{ p+1 }{ 2}  \right) \delta ^{ 2} s^{ -( p-1)} H_{2}
^{ p} + 2\delta  \varepsilon^{   p( p-1 ) } H_{3}^{ \frac{ p+1 }{ 2}
} + 2 \delta s H_{3}^{ \frac{ p+1 }{ 2} } + 2 H_{3}.  $$
 By  \eqref{3.8},  we have   $$  H_{ 3} \geq  C_{ 0}/4,   $$
   $$  \frac{ 1}{ 8} K_{ 0} s^{  -(p-1) } H_{2}^{ p} =  \frac{ 1}{ 8} K_{ 0} s H_{3}^{ p}
    \geq  \frac{ 1}{ 8} K_{ 0} s_{ 0} \left( \frac{ C_{ 0} }{ 4} \right )^{ p-1} H_{3}.   $$
As long as $ s_{ 0} $ sufficiently large, we have
    $  \displaystyle \frac{ 1}{ 8} K_{ 0} s_{ 0} \left( \displaystyle\frac{ C_{ 0} }{ 4} \right )^{ p-1} >1
    $,
 so we have
  $$  H_{ 3}<  \displaystyle \frac{ 1}{ 8} K_{ 0} s^{  -(p-1) } H_{2}^{ p},                  $$
 in the same way,   as long as $ \delta $ is sufficiently small,   $ \varepsilon \leq 1
 $,       we have
 $$  2\delta  \varepsilon^{   p( p-1 ) }  H_{3}^{ \frac{ p+1
}{ 2} } + 2 \delta s   H_{3}^{ \frac{ p+1 }{ 2} } \leq \displaystyle
\frac{ 1}{ 8} K_{ 0} s  H_{3}^{ p}=  \displaystyle \frac{ 1}{ 8} K_{
0} s^{ 1-p}  H_{2}^{ p} =  \displaystyle \frac{ 1}{ 8} K_{ 0} s^{ -(
p-1) }  H_{2}^{ p}( s ).            $$

  Similarly,  when  $ \varepsilon \leq 1  $,  and $ \delta $ is sufficiently
  small,  we have

  $$   \varepsilon^{   p( p-1 ) }  \left(  \frac{ p+1 }{ 2}  \right) \delta^{ 2}
  < \displaystyle \frac{ 1}{ 8} K_{ 0},   $$
  so we have
$$   \varepsilon^{   p( p-1 ) } H_{2}''( s ) +  2  H_{2}'( s )
    <   \displaystyle \frac{ 1}{ 2} K_{ 0} s^{  -(p-1) } H_{2}^{ p} (s ) <
     K_{ 0} s^{  -(p-1) } H_{2}^{ p} (s ),  $$
 $$     H_{2}( s_{ 0} )= C_{ 0}s_{ 0}/4,                      $$
  and  as long as $ \delta $ is sufficiently small,  we have
 $$     H_{2}'( s_{ 0} ) =  \displaystyle\frac{ C_{ 0} }{ 4}  +
 \delta s_{ 0}  \left( \displaystyle\frac{ C_{ 0} }{ 4} \right )^{ \frac{ p+1 }{ 2}
 } < C_{ 0} . $$
 So we obtain that    $  H_{2}( s ) $ satisfies that
$$
  \left\{
     \begin{array}{ll}
      \varepsilon^{   p( p-1 ) }  H_{2}''( s )+ 2  H_{2}'( s )  < K_{ 0} s^{  -(p-1) } H_{2}^{ p} (s ),
                \    \        \cr\noalign{\vskip4mm}
         H_{2}( s_{ 0} ) =  C_{ 0}s_{ 0}/4,   \    \       \       H_{2}'( s_{ 0} ) < C_{ 0}.
\end{array}
   \right.
$$
So by  Lemma \ref{thm:2.1}, we have
$$  H_{1}( s ) > H_{2}( s ),   \     \       s\geq s_{ 0},     $$
 and  since  $ H_{2}( s ) =s H_{3}( s )$,  and   $H_{3}( s ) $
 satisfies the Ricatti equation,  so there exists a positive
 constant $ s_{ 1} $ which is independent of  $ \varepsilon $,  such that
 $ H_{3}  $ blows up at time  $ s_{ 1} $, so   $ H_{2}  $ is infinity
at $ s_{ 1} $,  therefore $ H_{1}  $ definitely blow up before the
time $s_{ 1} $.  By the construction of $ H_{1}  $, we know that   $
H_{0} $ must blow up before the time  $ \varepsilon^{  -p( p-1 )
}s_{ 1} $, still by the construction of $ H_{0}  $,   we know that
$H(t)$ and then  the solution  $ u$ of \eqref{1.1}  must blow up
before the time  $  \exp ( \varepsilon^{  -p( p-1 ) }s_{ 1}) -2 $.
Thus we prove the upper bound of the lifespan in the theorem.

 The proof of Theorem 1.1 is complete.

\section*{Acknowledgments.}
The first author would like to thank Professor B. Yordanov for
helpful discussions.

This work is supported by the National Natural Science Foundation of
China (10728101),  the Basic Research Program of China  (No.
2007CB814800),  the Doctoral Foundation of the Ministry of Education
of China,  the `111' Project (B08018) and SGST 09DZ2272900.

\end{document}